\documentclass[11pt]{article}
\usepackage{amsfonts}
\usepackage{amsmath}
\usepackage{amssymb}
\usepackage{graphicx}
\usepackage[dvipsnames,usenames]{color}
\usepackage[pdfstartview=FitH ]{hyperref}
\usepackage[hyperpageref]{backref}
\begin{document}

\baselineskip 18pt
\def\o{\over}
\def\e{\varepsilon}
\title{\Large\bf A\ \ Note\ \ on\ \ Terence\ \ Tao's\ \ Paper\ \ ``On
\ \ the\\
Number\ \ of\ \ Solutions\ \ to\ \ $\bf {4\o p}={1\o n_1}
+{1\o n_2}+{1\o n_3}$''}
\author{Chaohua\ \  Jia}
\date{}
\maketitle {\small \noindent {\bf Abstract.} For the positive integer
$n$, let $f(n)$ denote the number of positive integer solutions
$(n_1,\,n_2,\,n_3)$ of the Diophantine equation
$$
{4\o n}={1\o n_1}+{1\o n_2}+{1\o n_3}.
$$
For the prime number $p$, $f(p)$ can be split into $f_1(p)+f_2(p),$
where $f_i(p)(i=1,\,2)$ counts those solutions with exactly $i$ of
denominators$n_1,\,n_2,\,n_3$ divisible by $p.$

Recently Terence Tao proved that
$$
\sum_{p< x}f_2(p)\ll x\log^2x\log\log x.
$$
with other results. But actually only the upper bound $x\log^2x\log\log^2x$ 
can be obtained in his discussion. In this note we shall use an elementary 
method to save a factor $\log\log x$ and recover the above estimate.
}

\vskip.5in
\noindent{\bf 1. Introduction}

For the positive integer $n$, let $f(n)$ denote the number of positive
integer solutions $(n_1,\,n_2,\,n_3)$ of the Diophantine equation
$$
{4\o n}={1\o n_1}+{1\o n_2}+{1\o n_3}.
$$
Erd\" os and Straus conjectured that for all $n\geq 2,\, f(n)>0.$ It
is still an open problem now although there are some partial results.

In 1970, R. C. Vaughan[2] showed that the number of $n< x$ for which
$f(n)=0$ is at most $x\exp(-c\log^{2\o 3}x),$ where $x$ is sufficiently
large and $c$ is a positive constant.

Recently Terence Tao[1] studied the situation in which $n$ is the prime
number $p.$ He gave lower bound and upper bound for the mean value of
$f(p).$ Precisely, he split $f(p)$ into $f_1(p)+f_2(p),$ where $f_i(p)
(i=1,\,2)$ counts those solutions with exactly $i$ of denominators
$n_1,\,n_2,\,n_3$ divisible by $p.$ He proved that
\begin{align}
x\log^2x\ll\sum_{p< x}f_2(p)\ll x\log^2x\log\log x
\end{align}
and
$$
x\log^2x\ll\sum_{p< x}f_1(p)\ll x\exp({c\log x\o \log\log x}),
$$
where $p$ denotes the prime number, $x$ is sufficiently large and $c$
is a positive constant. Then he conjectured that for $i=1,\,2$,
\begin{align}
\sum_{p< x}f_i(p)\ll x\log^2x.
\end{align}

But actually Terence Tao[1] only proved
\begin{align}
\sum_{p< x}f_2(p)\ll x\log^2x\log\log^2x,
\end{align}
since there was an error in his discussion. In this note we shall use
an elementary method to save a factor $\log\log x$ and recover the
upper bound in the right side of (1).

{\bf Theorem.} Let $p$ denote the prime number. Then for sufficiently
large $x$, we have
$$
\sum_{p< x}f_2(p)\ll x\log^2x\log\log x.
$$

\vskip.3in
\noindent{\bf 2. The proof of Theorem}

{\bf Lemma 1.} If $\varphi(n)$ is the Euler totient function, then
$$
\varphi(n)\gg {n\o g(n)}.
$$
Here
$$
g(n)=\prod_{p|n}(1+{1\o p})=\sum_{d|n}{\mu^2(d)\o d},
$$
where $\mu(d)$ is the M\"obius functions.

{\bf Proof.} We know that
$$
\varphi(n)=n\prod_{p|n}(1-{1\o p}).
$$
Then
$$
\varphi(n)=n\,{\prod_{p|n}(1-{1\o p^2})\o \prod_{p|n}(1+{1\o p})}
\geq {n\o g(n)}\prod_p(1-{1\o p^2})\gg {n\o g(n).}
$$
It is easy to see
$$
g(n)=\sum_{d|n}{\mu^2(d)\o d}.
$$

{\bf Lemma 2.} If $x\geq 1$, then
$$
\sum_{x< n\leq 2x}{1\o \varphi(n)}\ll 1.
$$

{\bf Proof.} By Lemma 1, we have
\begin{align*}
\sum_{x< n\leq 2x}{1\o \varphi(n)}&\ll\sum_{x< n\leq 2x}{g(n)\o n}\\
&=\sum_{x< n\leq 2x}{1\o n}\sum_{d|n}{\mu^2(d)\o d}\\
&=\sum_{d\leq 2x}{\mu^2(d)\o d}\sum_{\substack{x< n\leq 2x\\ d|n}}{1\o n}\\
&=\sum_{d\leq 2x}{\mu^2(d)\o d^2}\sum_{{x\o d}< l\leq {2x\o d}}{1\o l}\\
&\ll\sum_{d\leq 2x}{\mu^2(d)\o d^2}\ll 1.
\end{align*}

{\bf Lemma 3.} Let $p$ denote the prime number. Then the functions
$f_2(p)$ is equal to three times the number of triples $(a,\,b,\,c)$
of positive integers such that
$$
(a,\,b)=1,\ \ c|a+b,\ \ 4ab|p+c.
$$

One can see Proposition 1.2 of [1].

By some transformation, Terence Tao[1] got
\begin{align*}
\sum_{p< x}f_2(p)&\ll\sum_{1\leq i\leq {1\o 2}\log_2 x}\sum_{i\leq
j\leq \log_2 x-i}{x\o 1+\log_2 x-i-j}\cdot\\
&\cdot\sum_{2^i< a\leq 2^{i+1}}\sum_{\substack{2^j< b\leq 2^{j+1}
\\ (a,\, b)=1}}{d(a+b)\o \varphi(a)\varphi(b)}.
\end{align*}
Here $d(n)$ is the divisor function. It is necessary to keep the
condition $(a,\,b)=1.$

Now we consider the estimate for the sum
\begin{align}
\sum_{V< a\leq 2V}{1\o \varphi(a)}\sum_{\substack{W< b\leq 2W\\ (a,\,b)
=1}}{d(a+b)\o \varphi(b)},
\end{align}
where $1\leq V\leq W\leq x.$

Let
\begin{align}
S(a,\,W)=\sum_{\substack{W< b\leq 2W\\ (a,\,b)=1}}{d(a+b)\o \varphi(b)}.
\end{align}
Then Lemma 1 yields that
\begin{align*}
S(a,\,W)&\ll\sum_{\substack{W< b\leq 2W\\ (a,\,b)=1}}d(a+b)\cdot{g(b)\o b}\\
&\ll{1\o W}\sum_{\substack{W< b\leq 2W\\ (a,\,b)=1}}d(a+b)g(b)\\
&={1\o W}\sum_{\substack{W+a< k\leq 2W+a\\ (k,\,a)=1}}d(k)g(k-a)\\
&={1\o W}\sum_{\substack{W+a< rl\leq 2W+a\\ (rl,\,a)=1}}g(rl-a)\\
&\leq{2\o W}\sum_{\substack{r\leq\sqrt{2W+a}\\ (r,\,a)=1}}
\sum_{\substack{l\\ W+a< rl\leq 2W+a\\ (l,\,a)=1}}g(rl-a)\\
&\ll{1\o W}\sum_{\substack{r\leq\sqrt{2W+a}\\ (r,\,a)=1}}
\sum_{\substack{W< n\leq 2W\\ n\equiv -a({\rm mod}\,r)\\ (n,\,a)=1}}g(n)\\
&\leq{1\o W}\sum_{\substack{r\leq\sqrt{2W+a}\\ (r,\,a)=1}}
\sum_{\substack{W< n\leq 2W\\ n\equiv -a({\rm mod}\,r)}}g(n).
\end{align*}

Since $(r,\,a)=1,\,n\equiv-a({\rm mod}\,r)\Longrightarrow (n,\,r)=1.$
Then
\begin{align*}
\sum_{\substack{W< n\leq 2W\\ n\equiv -a({\rm mod}\,r)}}g(n)&=
\sum_{\substack{W< n\leq 2W\\ n\equiv -a({\rm mod}\,r)\\ (n,\,r)=1}}g(n)\\
&=\sum_{\substack{W< n\leq 2W\\ n\equiv -a({\rm mod}\,r)\\ (n,\,r)=1}}
\sum_{d|n}{\mu^2(d)\o d}\\
&=\sum_{\substack{d\leq 2W\\ (d,\,r)=1}}{\mu^2(d)\o d}
\sum_{\substack{W< n\leq 2W\\ n\equiv -a({\rm mod}\,r)\\ (n,\,r)=1\\
d|n}}1\\
&=\sum_{\substack{d\leq 2W\\ (d,\,r)=1}}{\mu^2(d)\o d}
\sum_{\substack{{W\o d}< k\leq {2W\o d}\\ (k,\,r)=1\\ dk\equiv -a({\rm mod}\,r)}}1\\
&=\sum_{\substack{d\leq 2W\\ (d,\,r)=1}}{\mu^2(d)\o d}
\sum_{\substack{{W\o d}< k\leq {2W\o d}\\ (k,\,r)=1\\ k\equiv -{\bar d}
a({\rm mod}\,r)}}1\\
&\leq\sum_{\substack{d\leq 2W\\ (d,\,r)=1}}{\mu^2(d)\o d}
\sum_{\substack{{W\o d}< k\leq {2W\o d}\\ k\equiv -{\bar d}
a({\rm mod}\,r)}}1,
\end{align*}
where ${\bar d}$ is an integer such that ${\bar d}d\equiv 1({\rm mod}\,r).$

We have
$$
\sum_{\substack{{W\o d}< k\leq {2W\o d}\\ k\equiv -{\bar d}
a({\rm mod}\,r)}}1
\ll{W\o dr}+1.
$$
Thus
\begin{align*}
\sum_{\substack{W< n\leq 2W\\ n\equiv -a({\rm mod}\,r)}}g(n)&\ll
\sum_{\substack{d\leq 2W\\ (d,\,r)=1}}{\mu^2(d)\o d}({W\o dr}+1)\\
&\leq {W\o r}\sum_{d\leq 2W}{\mu^2(d)\o d^2}+\sum_{d\leq 2W}{\mu^2(d)\o d}\\
&\ll {W\o r}+\log 2W.
\end{align*}
It follows that
\begin{align*}
S(a,\,W)&\ll {1\o W}\sum_{\substack{r\leq\sqrt{2W+a}\\ (r,\,a)=1}}({W\o r}+\log 2W)\\
&\leq{1\o W}\sum_{r\leq 2\sqrt{W}}({W\o r}+\log 2W)\\
&\ll\log 2W.
\end{align*}
By Lemma 2, we have
\begin{align*}
&\ \sum_{V< a\leq 2V}{1\o \varphi(a)}\sum_{\substack{W< b\leq 2W\\ (a,\,b)=1}}
{d(a+b)\o \varphi(b)}\\
&\ll\sum_{V< a\leq 2V}{1\o \varphi(a)}\cdot\log x\\
&\ll\log x.
\end{align*}
Therefore
\begin{align}
\sum_{p<x}f_2(p)\ll x\log x\sum_{1\leq i\leq{1\o 2}\log_2x}\sum_{i
\leq j\leq \log_2x-i}{1\o 1+\log_2x-i-j}.
\end{align}

We have
\begin{align*}
&\ \sum_{1\leq i\leq{1\o 2}\log_2x}\sum_{i\leq j\leq \log_2x-i}{1\o
1+\log_2x-i-j}\\
&\leq\sum_{1\leq i\leq{1\o 2}\log_2x}\sum_{1\leq h\leq \log_2x-2i+2}
{1\o h}\\
&\ll\sum_{1\leq i\leq{1\o 2}\log_2x}\log(\log_2x-2i+4)\\
&\ll\sum_{1\leq k\leq{1\o 2}\log_2x+1}\log(2k+8)\\
&\ll\log x\log\log x.
\end{align*}

So far the proof of Theorem is finished.

Similar discussion can yield
$$
\sum_{1\leq i\leq{1\o 2}\log_2x}\sum_{i\leq j\leq \log_2x-i}{1\o
1+\log_2x-i-j}\gg \log x\log\log x.
$$
In [1],
$$
\sum_{1\leq i\leq{1\o 2}\log_2x}\sum_{i\leq j\leq \log_2x-i}{1\o
1+\log_2x-i-j}\ll \log x
$$
is proved, where a factor $\log\log x$ is lost. From the above
discussion, it seems reasonable to conjecture
\begin{align}
x\log^2x\log\log x\ll\sum_{p< x}f_2(p).
\end{align}


\bigskip

\

Institute of  Mathematics, Academia Sinica, Beijing 100190, P. R.
China

E-mail: jiach@math.ac.cn
\end{document}